\documentclass[11pt, a4paper ,reqno]{amsart}
\usepackage{syntonly}
\usepackage{pstricks,subfigure,epsfig}
\usepackage{amsmath,amsfonts,amssymb,amsthm,euscript,upgreek}
\usepackage{enumerate}
\usepackage{multirow}
\usepackage{appendix}

\newcommand{\R}{\mathbb{R}}                     

\newcommand{\CP}{\mathbb{C}\mathrm{P}}

\newcommand{\CH}{\mathbb{C}\mathrm{H}}
\newcommand{\ol}{\mathrm{Hol}}
\newcommand{\hilb}{\mathcal{H}}
\newcommand{\f}{\rightarrow}                    
\newcommand{\C}{\mathbb{C}}                     
\newcommand{\de}{\partial}                      

\newcommand\aut{\mathrm{Aut}}
\newcommand\isom{\mathrm{Isom}}

\newcommand{\K}{K\"{a}hler}

\newtheorem{theor}{Theorem}

\newtheorem{lem}[theor]{Lemma}
\newtheorem{cor}[theor]{Corollary}

\newtheorem{ex}{Example}
\newtheorem{remar}[theor]{Remark}

\begin{document}
\title{Balanced metrics on Cartan and Cartan--Hartogs domains}
\author[A. Loi, M. Zedda]{Andrea Loi, Michela Zedda}
\address{Dipartimento di Matematica e Informatica, Universit\`{a} di Cagliari,
Via Ospedale 72, 09124 Cagliari, Italy}
\email{loi@unica.it; michela.zedda@gmail.com  }
\thanks{
The first author was supported  by the M.I.U.R. Project \lq\lq Geometric
Properties of Real and Complex Manifolds'';
the second author was  supported by  RAS
through a grant financed with the ``Sardinia PO FSE 2007-2013'' funds and 
provided according to the L.R. $7/2007$.}
\date{}

\begin{abstract}
This paper consists of two results dealing with balanced metrics (in S. Donaldson terminology) on nonconpact complex  manifolds.
In the first one we describe all balanced metrics on Cartan domains. In the second one we show that the only Cartan--Hartogs domain 
which admits a balanced metric is the complex hyperbolic space. By combining these results with those obtained in \cite{articwall}
we also provide the first  example of complete,   \K\--Einstein  and projectively induced metric $g$  such that $\alpha g$ is not balanced for all $\alpha >0$.  
\end{abstract}

\maketitle

\section{Introduction}
Let $\Omega\subset \C^d$ be a Cartan domain, i.e. an irreducible bounded symmetric domain, of complex dimension $d$ and genus $\gamma$. For all positive real numbers $\mu$ consider the family of \emph{Cartan-Hartogs} domains 
\begin{equation}\label{defm}
M_{\Omega}(\mu)=\left\{(z,w)\in \Omega\times\C,\ |w|^2<N_\Omega^\mu(z,z)\right\},
\end{equation}
where $N_\Omega(z,z)$ is the  {\em generic norm} of $\Omega$, i.e.
\begin{equation}\label{genericnorm}
N_{\Omega}(z, z)=(V(\Omega)K(z, z))^{-\frac{1}{\gamma}},
\end{equation}
where $V(\Omega)$ is the total volume of $\Omega$ with respect to the Euclidean measure of the ambient complex Euclidean space and $K(z, z)$ is its Bergman kernel.

The domain $\Omega$ is called  the {\em  base} of the Cartan--Hartogs domain 
$M_{\Omega}(\mu)$ (one also  says that 
$M_{\Omega}(\mu)$  is based on $\Omega$).
Consider on $M_{\Omega}(\mu)$ the metric $g(\mu)$  whose associated K\"ahler form $\omega(\mu)$ can be described by the (globally defined)
K\"ahler potential centered at the origin
\begin{equation}\label{diastM}
\Phi(z,w)=-\log(N_{\Omega}^\mu(z,z)-|w|^2).
\end{equation}
These domains have been considered by several authors (see e.g. \cite{roos} and references therein).
In \cite{articwall} the authors of the present paper  study when $\left(M_{\Omega}(\mu),\alpha\,g(\mu)\right)$, for a positive constant $\alpha$, admits a holomorphic and isometric (from now on \emph{\K}) immersion $f$ into the infinite dimensional complex projective space $\CP^\infty$, i.e. $f^*g_{FS}=\alpha\,g(\mu)$, where $g_{FS}$ denotes the Fubini--Study metric on $\CP^\infty$ (when such a K\"ahler immersion exists, we say also that the metric is \emph{projectively induced}). Recall that given homogeneous coordinates $[Z_0,\dots, Z_j,\dots]$ on $\CP^\infty$, $g_{FS}$ is the K\"ahler metric whose associated K\"ahler form $\omega_{FS}$ can be described in the open set $U_0=\{Z_0\neq 0\}$ by
$\omega_{FS}=\frac{i}{2}\de\bar \de\Phi_{FS}$, where $\Phi_{FS}= \log(1+\sum_{j=1}^{\infty}|z_j|^2)$  for $z_j=\frac{Z_j}{Z_0}$, $j=1,\dots$, affine coordinates on $U_0$.
The  main results obtained in  \cite{articwall} can be summarized in the following theorem (see also next section for a more detailed description of the Wallach set
$W(\Omega)$ and for the definition of the integer  $a$ appearing in (c)).

\vskip 0.3cm

\noindent
{\bf Theorem LZ}\label{wallach}
{\em Let $\Omega\subset \C^d$ be a Cartan domain of rank $r$, genus $\gamma$ and dimension $d$ and let $g_B$  be its Bergman metric. 
Then the  following results hold true:
 \begin{itemize}
 \item [(a)]
 $(\Omega,\beta  g_B)$, $\beta >0$, admits a equivariant  K\"ahler immersion into $\CP^\infty$ if and only if $\beta \gamma$ belongs to $W(\Omega)\setminus \{0\}$;
 \item [(b)]
 the metric  $\alpha g(\mu)$, $\alpha >0$, on the Cartan--Hartogs domain $M_{\Omega}(\mu)$
is projectively induced if and only if  $(\alpha +m)\frac{\mu}{\gamma}g_B$ is 
projectively induced for every integer $m\geq 0$; 
\item[(c)]
Let   $\mu_0=\gamma/(d+1)$  and $\Omega\neq\CH^d$. Then the metric $\alpha g(\mu_0)$ on   $M_\Omega(\mu_0)$ is  \K\--Einstein,   complete, nonhomogeneous and 
projectively induced  for all positive real number $\alpha\geq \frac{(r-1)(d+1)a}{2\gamma}$.
\end{itemize}}

In this paper we study balanced metrics (in S. Donaldson's terminology) on Cartan and Cartan--Hartogs domains.
The main results are  the following two  theorems. In the first one we describe all balanced metrics on Cartan's domains, while  the second one  can be viewed as a characterization of the complex hyperbolic space among Cartan--Hartogs domains, in terms of balanced metrics
(cfr. Example \ref{exhyp} below).

\begin{theor}\label{cartanbalanced}
Let $\Omega$ be a Cartan domain of genus $\gamma$ equipped with its Bergman metric $g_B$.
The metric $\beta g_B$, $\beta >0$,   is balanced if and only if $\beta>\frac{\gamma-1}{\gamma}$. 
\end{theor}

\begin{theor}\label{cartanhartogs}
Let $M_\Omega(\mu)$ be a Cartan-Hartogs domain based on the Cartan domain $\Omega\subset \C^d$. The metric $\alpha g(\mu)$ on $M_\Omega(\mu)$ is balanced if and only if $\alpha>d+1$ and $M_\Omega(\mu)$ is holomorphically isometric to the complex hyperbolic space $\CH^{d+1}$, namely $\Omega =\CH^d$ and $\mu=1$. 
\end{theor}

 By combining these results with  (c) in Theorem LZ  
we also obtain the first  example of {\em complete,   \K\--Einstein  and projectively induced metric $g$  such that $\alpha g$ is not balanced for
$\alpha$ varying in a continuous subset of the real numbers}. This is   expressed by the following corollary.

\begin{cor}
Let $\Omega\subset \C^d$ be a Cartan domain of genus $\gamma$ equipped with its Bergman metric $g_B$.
Let   $\mu_0=\gamma/(d+1)$  and $\Omega\neq\CH^d$. Then the metric $\alpha g(\mu_0)$ on   $M_\Omega(\mu_0)$ is  complete,   \K\--Einstein 
projectively induced and not balanced   for all $\alpha\geq \frac{(r-1)(d+1)a}{2\gamma}$.
\end{cor}

The paper consists in other three  sections. In Section \ref{balancedmetrics} we recall the definition  of   balanced metrics.
 In  Section \ref{balancedcartan}  we describe all
balanced metrics  on Cartan domains and prove Theorem \ref{cartanbalanced}. Finally Section \ref{balancedcartanhartogs} 
is dedicated to the proof of  Theorem \ref{cartanhartogs}.

\section{Balanced metrics}\label{balancedmetrics}
Let $M$ be a $n$-dimensional complex manifold  endowed with a K\"ahler metric $g$
and let $\omega$ be the \K\ form associated to $g$, i.e. $\omega (\cdot ,\cdot  )=g(J\cdot, \cdot)$.
Assume that  the metric $g$ can be described by a
strictly plurisubharmonic real valued function $\Phi :M\rightarrow\R$, called a {\em K\"ahler potential} for $g$, 
i.e.
$\omega=\frac{i}{2}\de\bar\de \Phi$.
 
Let $\hilb_\Phi$ be the weighted Hilbert space of square integrable holomorphic functions on $(M, g)$, with weight $e^{-\Phi}$, namely
\begin{equation}\label{hilbertspace}
\hilb_\Phi=\left\{ f\in\ol(M) \ | \ \, \int_M e^{-\Phi}|f|^2\frac{\omega^n}{n!}<\infty\right\},
\end{equation}
where $\frac{\omega^n}{n!}=\det(\de\bar \de \Phi)\frac{\omega_0^n}{n!}$ is the volume form associated to $\omega$ and $\omega_0=\frac{i}{2}\sum_{j=1}^{n} dz_j\wedge d\bar z_j$ is the standard K\"ahler form on $\C^n$.
If $\hilb_\Phi\neq \{0\}$ we can pick an orthonormal basis $\{f_j\}$ and define its reproducing kernel by
$$K_{\Phi}(z, z)=\sum_{j=0}^N |f_j(z)|^2, $$
where $N+1$ denotes the complex dimension of $\hilb_\Phi\neq \{0\}$.
Consider the function 
\begin{equation}\label{epsilon}
\varepsilon_g(z)=e^{-\Phi(z)}K_{\Phi}(z, z). 
\end{equation} 
As suggested by the notation it is not difficult to verify  that  $\varepsilon _g$ depends only on the metric $g$ and not on the choice of the K\"ahler potential $\Phi$ 
(which is defined up to an addition with the real part of a holomorphic function on $M$) 
or on  the orthonormal basis chosen.
\vskip 0,3cm

\noindent {\bf Definition.} The metric $g$ is \emph{balanced} if  the function $\varepsilon_g$ is a positive  constant.
\vskip 0,3cm

A balanced metric $g$ on $M$ can be viewed as  a particular  projectively induced  \K\ metric for which the K\"ahler immersion $f\!: M\f \CP^N, N\leq\infty,\ x\mapsto [s_0(x), \dots ,s_j(x), \dots]$, is given by the  orthonormal basis $\{f_j\}$ of the Hilbert space $\hilb_{\Phi}$. Indeed the map $f$ is well-defined since $\varepsilon_g$ is a positive constant
and hence for all $x\in M$ there exists $\varphi\in \hilb_{\Phi}$ such that $\varphi(x)\neq 0$. Moreover, 
\begin{equation}
\begin{split}
f^*\omega_{FS}=&\frac{i}{2}\de\bar\de\log\sum_{j=0}^\infty|f_j(z)|^2\\
=&\frac{i}{2}\de\bar\de\log K_{ \Phi} (z, z)\\
=&\frac{i}{2}\de\bar\de\log \varepsilon_{\,g}+\frac{i}{2}\de\bar\de\log e^{\Phi}\\
=&\frac{i}{2}\de\bar\de\log \varepsilon_{g}+   \,\omega.\nonumber
\end{split}
\end{equation}
Hence if $g$ is balanced the map $f$ is isometric.

In the literature the function $\varepsilon_g$ was first introduced under the name of $\eta$-{\em function} by J. Rawnsley in \cite{rawnsley}, later renamed as $\varepsilon$-{\em function} in \cite{CGR}. The map $f$ is called in \cite{CGR} the  {\em coherent states map}.  It  plays a fundamental role in the geometric quantization and quantization by deformation  of a K\"ahler manifold. It also related to the 
Tian-Yau-Zelditch asymptotic expansion (see   \cite{graloi},  \cite{quant}, \cite{constscal} and references therein).

\begin{ex}\rm\label{exhyp}
 Notice that a projectively induced metric is not always balanced. For example, in \cite{ca} E. Calabi shows that the complex hyperbolic space $(\CH^d,  \alpha \,g_{hyp})$,
 endowed with a positive multiple of the hyperbolic metric $g_{hyp}$, is projectively induced for all $\alpha>0$. (Here  $\CH^d=\{z\in \C^d\ | \ |z|^2<1\}$
 and the \K\ form $\omega_{hyp}$ associated to  $g_{hyp}$ is given by  
 $\omega_{hyp}=-\frac{i}{2}\partial\bar\partial\log (1-|z|^2)$).
  Althought, it is well-known that the weighted Hilbert space of square integrable holomorphic functions on $(\CH^d,  \alpha \,g_{hyp})$, 
i.e.
$$\hilb_{\alpha\Phi_{hyp}}=\left\{ \varphi\in\ol(\CH^d), \int_{\CH^d}\left(1-|z|^2\right)^{ \alpha -(d+1)}|\varphi|^2\frac{\omega_0^d}{d!}<\infty \right\},$$
is equal to  $\{0\}$ for all $\alpha\leq d$. Similar considerations can be done for all Cartan domains (see Remark \ref{remarbalsym} below).
\end{ex}

\begin{remar}\rm
 The definition of balanced metrics
 was originally given by S. Donaldson \cite{donaldson} in the case of  a compact polarized \K\ manifold $(M, g)$ and generalized in \cite{arezzoloi}
 (see also \cite{cucculoibal}, \cite{englisweigh}, \cite{grecoloi}) to the noncompact case. Here we give only the definition for those \K\ metrics which admits a globally defined potential as the Cartan and Cartan--Hartogs domains treated in this paper.
\end{remar}

\section{Balanced metrics on Cartan domains}\label{balancedcartan}
Let $(\Omega, \beta g_B)$, $\beta >0$,  denote a Cartan domain, i.e. an irreducible bounded symmetric domain of $\C^d$ endowed with a positive multiple of its Bergman metric $g_B$. Recall that $g_B$ is the K\"ahler metric on $\Omega$ whose associated K\"ahler form $\omega_B$ is given by $\omega_B=\frac{i}{2}\de\bar\de\log K(z, z)$, where  $K(z, z)$ is the reproducing kernel for the Hilbert space
$$\hilb =\left\{\varphi\in\ol(\Omega),\ \int_\Omega |\varphi|^{2d}\ \frac{\omega_0^d}{d!}<\infty\right\},$$
where $\omega_0=\frac{i}{2} \sum_{j=1}^d dz_j\wedge d\bar z_j$ is the standard K\"ahler form of $\C^d$.
A bounded symmetric domain $(\Omega,\alpha g_B)$ is uniquely determined by a triple of integers $(r,a,b)$, where $r$ represents the rank of $\Omega$ and $a$ and $b$ are positive integers. The genus $\gamma$ of $\Omega$ is defined by  $\gamma=(r-1)a+b+2$. 
The table below summarizes the numerical invariants and the dimension of $\Omega$ according to its type (for a more detailed description of this invariants, which is not necessary in our approach, see e.g. \cite{arazy}).
\begin{table}[!h]\label{inv}
\caption{Bounded symmetric domains, invariants and dimension.}
\vskip -0.5cm
\centering      
\begin{tabular}{cccccc} \label{inv}
$\textrm{}$&&&&&\\
\hline                       
 Type& $r$ & $a$ & $b$ & $\gamma$  & dimension\\
\hline                    
$\Omega_1[m,n]$ &$m$& $2$&$n-m$&$n+m $&$nm$  \\
$\Omega_2[n]$&  $n$& $1$&$0$& $n+1$ &$n(n+1)/2$ \\  \vspace{-0.2cm}
\multirow{2}{*}{$\Omega_3[n]$} & \multirow{2}{*}{$[n/2]$}&\multirow{2}{*}{$4$}& \small{$0$ ($n$ even)} &\multirow{2}{*}{$n-1$}&\multirow{2}{*}{$n(n-1)/2$} \\
&&&\small{$2$ ($n$ odd)}&&\\
$\Omega_4[n]$ & $2$&$n-2$& $0$&$n$&$n$\\
$\Omega_{V}[16]$&$2$&$6$&$4$&$12$&$16$\\
$\Omega_{VI}[27]$&$3$&$8$&$0$&$18$&$27$\\
 \hline    
\end{tabular} 
\end{table} 

We give now the definition of the  {\em Wallach set}  of a Cartan domain  $\Omega$, referring  the reader to  \cite{arazy}, \cite{faraut} and \cite{upmeier} for more details and results. The Wallach set, denoted by $W(\Omega)$, consists of all $\eta\in\C$ such that there exists a Hilbert space ${\hilb}_\eta$ whose reproducing kernel is   $K^{\frac{\eta}{\gamma}}$. This is equivalent to the requirement that $K^{\frac{\eta}{\gamma}}$ is positive definite, i.e.  for all $n$-tuples of  points $x_1,\dots,x_n$ belonging to $\Omega$ the $n\times n$ matrix $(K(x_{\alpha},x_{\beta})^{\frac{\eta}{\gamma}})$, is positive  {\em semidefinite}. It turns out (cfr. \cite[Cor. 4.4, p. 27]{arazy} and references therein) that $W(\Omega)$ consists only of real numbers and depends on two of the domain's invariants, $a$ and  $r$. More precisely we have
\begin{equation}\label{wallachset}
W(\Omega)=\left\{0,\,\frac{a}{2},\,2\frac{a}{2},\,\dots,\,(r-1)\frac{a}{2}\right\}\cup \left((r-1)\frac{a}{2},\,\infty\right).
\end{equation}
The set $W_{dis}=\left\{0,\,\frac{a}{2},\,2\frac{a}{2},\,\dots,\,(r-1)\frac{a}{2}\right\}$ and the interval $W_c= \left((r-1)\frac{a}{2},\,\infty\right)$
are called respectively  the {\em discrete} and {\em continuous} part   of the Wallach set of the domain 
$\Omega$.
\begin{remar}\label{rchimm}\rm
If $\Omega$ has rank $r=1$, namely $\Omega$ is the complex hyperbolic space $\CH^d$, then $g_B=(d+1)g_{hyp}$. In this case (and only in this case) $W_{dis}=\{0\}$ and $W_c=(0, \infty)$ (cfr. Example \ref{exhyp}).
\end{remar}

We can now proof Theorem \ref{cartanbalanced}.

\begin{proof}[Proof of Theorem \ref{cartanbalanced}]
Let $d$ denote the complex dimension of $\Omega$.
It follows by standard results  on bounded symmetric domains  (see e.g. \cite{faraut}) 
that the Hilbert space  
$${\hilb}_\beta =\left\{ \varphi\in\ol(\Omega) \ | \ \, \int_{\Omega}\frac{1}{K^{\beta}}|\varphi|^2\frac{\omega_B^{d}}{d!}<\infty\right\},$$
does not reduce to the zero dimensional space
 iff $\beta>\frac{\gamma-1}{\gamma}$.
 
Hence, in order to prove that $\beta\, g_B$ is balanced for $\beta>\frac{\gamma-1}{\gamma}$,  it  remains   to show that  for $\beta>\frac{\gamma-1}{\gamma}$ the map $h_{\beta}:\Omega\rightarrow \C P^{\infty}, \ x\mapsto [\dots , h_{\beta}^j(x), \dots]$, where  $\{h_{\beta}^j\}$ is an orthonormal basis  of ${\hilb}_\beta$, is a well-defined map  of $\Omega$ into $\CP^{\infty}$ and it is \K\,  i.e. 
\begin{equation}
h_{\beta}^*g_{FS}=\beta g_B. \nonumber
\end{equation}
To prove that $h_{\beta}$ is well-defined one needs to verify that for all $x\in\Omega$ there exists $\varphi\in{\hilb}_\beta$ such that $\varphi(x)\neq 0$. Assume, by contradiction, that there exists $x_0\in\Omega$ such that $\varphi(x_0)=0$ for all  $\varphi\in\hilb_\beta$. Write $\Omega=G/K$, where $G$ is a subgroup of $\aut(\Omega)\cap\isom(\Omega)$ which acts transitively on $\Omega$. Then for all $g\in G$, $\varphi\circ g$ is an element of ${\hilb}_\beta$ which, by assumption, vanishes on $x_0$. Thus $0= \varphi\circ g(x_0)= \varphi(gx_0)$ and since this holds true for all $g\in G$, $h_{\beta}$ is the zero function.  Hence $\hilb_\beta=\{0\}$, which is in contrast with the fatct that $\hilb_{\beta}\neq\{0\}$ for $\beta>\frac{\gamma -1}{\gamma}$.
In order to prove that $h_{\beta}$ is \K\ notice that  the function $\frac{\sum_{j=0}^{\infty}|h_{\beta}^j|^2}{K^{\beta}}$ is invariant by the group $G$ and hence constant. Hence
$$h_{\beta}^*\omega_{FS}=
\frac{i}{2}\partial\bar\partial\log \sum_{j=0}^{\infty}|h_{\beta}^j|^2= \beta\omega_B+\frac{i}{2}\partial\bar\partial\log \frac{\sum_{j=0}^{\infty}|h_{\beta}^j|^2}{K^{\beta}}=\beta\omega_B,$$
and we are done.
 \end{proof}
 \begin{remar}\rm\label{remarbalsym}
By Theorem \ref{cartanbalanced} the subset of the positive real numbers  $\beta$ for which $\beta \, g_B$ is balanced, i.e. $\left(\frac{\gamma-1}{\gamma},\infty\right)$,  is a proper subset of the continuous part $\left((r-1)\frac{a}{2},\infty\right)$ of the Wallach set $W(\Omega)$. Combining this remark  with (a) of Theorem LZ in the introduction one gets that for every Cartan domain there exists an infinite interval  of positive real numbers $\beta$ such that $\beta g_B$ is projectively induced  but not balanced. 
 \end{remar}
 \begin{remar}\rm\label{convergencetheorem1}
Observe that it follows by Theorem \ref{cartanbalanced} that, for all $\beta>\frac{\gamma-1}{\gamma}$, we have for some constant $\xi$
\begin{equation}\label{integralomega}
\int_{\Omega}N^{\gamma(\beta-1)}_\Omega\, h^j_{\beta}\,\bar h^k_{\beta}\,\omega_0^d=\xi\,\delta_{j,k},
\end{equation}
where $N_\Omega$ is the generic norm of $\Omega$ defined in (\ref{genericnorm}) and $h_\beta$ is the K\"ahler map defined in the proof of Theorem \ref{cartanbalanced}. In particular, the integral (\ref{integralomega}) is convergent and does not depend on $j$, $k$.
\end{remar}

\section{Balanced metrics on Cartan--Hartogs domains}\label{balancedcartanhartogs}
In order to prove Theorem \ref{cartanhartogs} we need the following two lemmata.
The first one gives an explicit description of the  \K\ immersions 
of a $d+1$-dimensional Cartan--Hartogs domain $(M_\Omega(\mu),\alpha g(\mu))$ into $\CP^{\infty}$ while the second one 
describes a necessary condition for the metric $\alpha g(\mu)$ to be balanced.

\begin{lem}\label{immersion}
If $f\!:M_\Omega(\mu)\f \CP^\infty$ is a holomorphic map such that $f^*\omega_{FS}=\alpha\,\omega(\mu)$ then up to unitary transformation of $\CP^\infty$
it is given by 
\begin{equation}\label{immf}
f=\left[ 1, s, h_{\frac{\mu\, \alpha}{\gamma}},\dots,\sqrt{\frac{(m+ \alpha-1)!}{(\alpha-1)!m!}}h_{\frac{\mu(\alpha +m)}{\gamma}}w^m,\dots\right],
\end{equation}
where $s=(s_1,\dots, s_m,\dots)$ with 
$$s_m=\sqrt{\frac{(m+ \alpha-1)!}{(\alpha-1)!m!}}w^m,$$
and for all $k>0$, $h_k=(h_k^1,\dots,h_k^j,\dots)$ is the sequence of holomorphic maps on $\Omega$ such that the immersion $\tilde h_k=(1,h_k^1,\dots, h_k^j,\dots)$, $\tilde h_k\!:\Omega\f\CP^\infty$, satisfies $\tilde h_k^*\omega_{FS}=k \omega_B$, i.e. 
\begin{equation}\label{ie}
1+\sum_{j=1}^{\infty}|h_k^j|^2=\frac{1}{N^{\gamma\, k}}.
\end{equation} 
\end{lem}
\begin{proof}
Since the immersion is isometric, by (\ref{diastM}) we have $f^*\Phi_{FS}=-\alpha\log(N_{\Omega}^\mu(z,z)-|w|^2)$, which is equivalent to 
$$\frac{1}{(N^{\mu}-|w|^2)^\alpha}=\sum_{j=0}^\infty |f_j|^2,$$
for $f=[f_0,\dots, f_j,\dots]$.
If we consider the power expansion around the origin of the left hand side with respect to $w$, $\bar w$, we get
\begin{equation}
\begin{split}
\sum_{k=1}^\infty \left[\frac{\de^{2k}}{\de w^k \de \bar w^k}\frac{1}{(N^{\mu}-|w|^2)^\alpha}\right]_{0}\frac{|w|^{2k}}{k!^2}=& \sum_{k=1}^\infty \left[\frac{\de^{2k}}{\de w^k \de \bar w^k}\frac{1}{(1-|w|^2)^\alpha}\right]_{0}\frac{|w|^{2k}}{k!^2}\\
=&\left(\sum_{k=0}^\infty |w|^2\right)^\alpha-1.\nonumber
\end{split}
\end{equation}
The power expansion with respect to $z$ and $\bar z$ reads
\begin{equation}
\begin{split}
\sum_{j,k}\left[\frac{\de^{|m_j|+|m_k|}}{\de z^{m_j} \de \bar z^{m_k}}\frac{1}{(N^{\mu}-|w|^2)^\alpha}\right]_{0}\frac{z^{m_j}\bar z^{m_k}}{m_j!m_k!}=&\sum_{j,k}\left[\frac{\de^{|m_j|+|m_k|}}{\de z^{m_j} \de \bar z^{m_k}}\frac{1}{N^{\mu\alpha}}\right]_{0}\frac{z^{m_j}\bar z^{m_k}}{m_j!m_k!}\\
=&\sum_{j=1}^\infty |h_{\frac{\mu\alpha}{\gamma}}^j|^2,\nonumber
\end{split}
\end{equation}
where the last equality holds since by (\ref{ie}) $\sum_{j=1}^\infty h_{\frac{\mu\alpha}{\gamma}}^j$ is the power expansion of $\frac{1}{N^{\mu\alpha}}-1$. Here we are using Calabi's multi index notation, namely we arrange every $d$-tuple of nonnegative integers as the sequence  $m_j=(m_{j,1},\dots,m_{j,d})$ with nondecreasing order, that is $m_0=(0,\dots,0)$, $|m_j|\leq |m_{j+1}|$, with $|m_j|=\sum_{\alpha=1}^d m_{j,\alpha}$. Further
$z^{m_j}$ denotes  the monomial in $d$ variables $\prod_{\alpha=1}^d z_\alpha^{m_{j,\alpha}}$ and $m_j!=m_{j,1}!\cdots m_{j,d}!$.

Finally, the power expansion with respect to $z$, $\bar z$, $w$, $\bar w$ reads
\begin{equation}
\begin{split}
&\sum_{m=1}^\infty\sum_{j,k}\left[\frac{\de^{|m_j|+|m_k|}}{\de z^{m_j} \de \bar z^{m_k}}\frac{\de^{2m}}{\de w^m \de \bar w^m}\frac{1}{(N^{\mu}-|w|^2)^\alpha}\right]_{0}\frac{z^{m_j}\bar z^{m_k}w^m\bar w^m}{m_j!m_k!m!^2}\\
=&\sum_{m=1}^\infty\sum_{j,k}\left[\frac{\de^{|m_j|+|m_k|}}{\de z^{m_j} \de \bar z^{m_k}}\frac{(m+\alpha-1)!}{(\alpha-1)!m!N^{\mu(\alpha+m)}}\right]_{0}\frac{z^{m_j}\bar z^{m_k}}{m_j!m_k!}\\
=&\sum_{m=1}^\infty\sum_{j=1}^\infty \frac{(m+\alpha-1)!}{(\alpha-1)!m!}|w|^{2m}|h_{\frac{\mu(\alpha+m)}{\gamma}}^j|^2,\nonumber
\end{split}
\end{equation}
where we are using (\ref{ie}) again. It follows by the previous power series expansions, that the map $f$ given by (\ref{immf}) is a \K\ immersion  of $(M_{\Omega}(\mu), \alpha g(\mu))$ into $\CP^\infty$. By Calabi's rigidity Theorem (cfr. \cite{ca}) all other \K\ immersions are given by $U\circ f$, where $U$ is a unitary transformation of $\C P^{\infty}$.
 \end{proof}

\begin{lem}\label{neccond}
If  $\alpha\,g(\mu)$ is  balanced then  $\alpha>d+1$ and $\alpha\mu>\gamma-1$. 
\end{lem}
\begin{proof}
Assume that $\alpha\, g(\mu)$ is balanced. Then it is projectively induced and by Lemma \ref{immersion}, up to unitary transformation of $\CP^\infty$, the K\"ahler immersion $f\!:M_\Omega(\mu)\f \CP^\infty$, $f=[f_0, \dots , f_j, \dots]$,  is given by (\ref{immf}). By Section \ref{balancedmetrics} $\{f_j\}_{j=0,1,\dots}$ is an orthonormal basis for the weighted Hilbert space
\begin{equation}\label{hilbertspaceexplicit}
\hilb_{\alpha\Phi}=\left\{ \varphi\in\ol(M_{\Omega}(\mu)) \ | \ \, \int_{M(\mu)}\left(N^\mu_\Omega-|w|^2\right)^\alpha|\varphi|^2\frac{\omega( \mu )^{d+1}}{(d+1)!}<\infty\right\},
\end{equation}
where up to the multiplication with a positive constant
 $$\frac{\omega(\mu)^{d+1}}{(d+1)!}=\frac{N_{\Omega}^{\mu(d+1)-\gamma}}{(N^\mu_\Omega-|w|^2)^{d+2}} \frac{\omega_0^{d+1}}{(d+1)!}.$$
Thus, in particular we have
\begin{equation}
\begin{split}
&\int_{M_{\Omega}(\mu)}(N^\mu_\Omega-|w|^2)^ \alpha f_j\bar f_k\frac{\omega(\mu)^{d+1}}{(d+1)!}=\\
&\int_{M_{\Omega}(\mu)}(N^\mu_\Omega-|w|^2)^{\alpha-(d+2)}N_{\Omega}^{\mu(d+1)-\gamma}f_j \bar f_k\frac{\omega_0^{d+1}}{(d+1)!}=\lambda\, \delta_{jk},\nonumber
\end{split}
\end{equation}
for some constant $\lambda$ indepentent from $j$ and $k$. It follows by (\ref{immf}) that the following integral 
\begin{equation}
\int_{M_\Omega(\mu)}(N^\mu_\Omega-|w|^2)^{\alpha-(d+2)}N^{\mu(d+1)-\gamma}|h^j_{\mu\alpha}|^2\frac{\omega_0^{d+1}}{(d+1)!},\nonumber
\end{equation}
is convergent. Passing to polar coordinates gets
\begin{equation}\label{integralm0}
\frac{\pi}{(d+1)!} \int_{\Omega}N_\Omega^{\mu(d+1)-\gamma} |h^j_{\mu\alpha}|^2 \int_0^{N_\Omega^\mu}(N^\mu_\Omega-\rho)^{\alpha-(d+2)}d\rho\,\omega_0^d.\nonumber
\end{equation}
The integral 
$$\int_0^{N_\Omega^\mu}(N^\mu_\Omega-\rho)^{\alpha-(d+2)}d\rho,$$ is convergent iff $\alpha-(d+2)>-1$, i.e. iff $\alpha>d+1$. Further, being $\alpha> d+1$, going on with computations gives 
\begin{equation}
\frac{\pi}{(d+1)!}\frac{1}{(\alpha-(d+2)+1)}\int_{\Omega}N^{\mu\alpha-\gamma}_\Omega |h^j_{\frac{\mu\alpha}{\gamma}}|^2\omega_0^d.\nonumber
\end{equation}
By Remark \ref{convergencetheorem1} this last integral converges and does not depends on $j$ iff $\alpha\mu>\gamma-1$, and we are done. 
\end{proof}
We are now in the position of proving Theorem \ref{cartanhartogs}.
\begin{proof}[Proof of Theorem \ref{cartanhartogs}]
Since by Theorem \ref{cartanbalanced} the hyperbolic metric $\alpha g_{hyp}$ is balanced iff $\alpha>d+1$, the sufficient condition is verified (recall that for the hyperbolic metric we have $\mu=1$ and $\gamma=d+2$).  For the necessary part, assume that $\alpha\, g(\mu)$ is balanced. By Lemma \ref{neccond}  we can assume $\alpha>d+1$ and $\alpha\mu>\gamma-1$.  Following the same approach as in Lemma \ref{neccond}, this gives that the integral
\begin{equation}
\int_{M_\Omega(\mu)}(N^\mu_\Omega-|w|^2)^{\alpha-(d+2)}N^{\mu(d+1)-\gamma}f_j\bar f_k\frac{\omega_0^{d+1}}{(d+1)!}\nonumber
\end{equation}
is zero for $j\neq k$ and does not depend on $j$ otherwise.
By (\ref{immf}) this implies that the following integral
\begin{equation}\label{integral}
\begin{split}
&\int_{M_\Omega(\mu)}(N^\mu_\Omega-|w|^2)^{\alpha-(d+2)}N^{\mu(d+1)-\gamma}\frac{(m+ \alpha-1)!}{(\alpha-1)!m!}|h^j_{\frac{\mu(\alpha +m)}{\gamma}}|^2|w|^{2m}\frac{\omega_0^{d+1}}{(d+1)!}=\\
&\frac{\pi}{(d+1)!} \frac{(m+ \alpha-1)!}{(\alpha-1)!m!}\int_{\Omega}N_\Omega^{\mu(d+1)-\gamma} |h^j_{\frac{\mu(\alpha +m)}{\gamma}}|^2 \int_0^{N_\Omega^\mu}(N^\mu_\Omega-\rho)^{\alpha-(d+2)}\rho^{m}\omega_0^d=\\
&\frac{\pi\, m!}{(d+1)!} \frac{(m+ \alpha-1)!}{(\alpha-1)!\,m!}\frac{1}{(\alpha-(d+2)+1)\cdots(\alpha-(d+2)+m)}\cdot\\
&\textrm{ }\qquad\qquad\qquad\qquad \cdot\int_{\Omega}N_\Omega^{\mu(d+1)-\gamma} |h^j_{\frac{\mu(\alpha +m)}{\gamma}}|^2 \int_0^{N_\Omega^\mu}(N^\mu_\Omega-\rho)^{\alpha-(d+2)+m}\omega_0^d=\\
&\frac{\pi}{(d+1)!} \frac{(m+ \alpha-1)!}{(\alpha-1)!}\frac{1}{(\alpha-(d+2)+1)\cdots(\alpha-(d+2)+m+1)}\cdot\\
&\textrm{ }\qquad\qquad\qquad\qquad\qquad\qquad\qquad\qquad\qquad \cdot\int_{\Omega}N^{\mu(\alpha +m)-\gamma}_\Omega |h^j_{\frac{\mu(\alpha +m)}{\gamma}}|^2\omega_0^d,
\end{split}
\end{equation}
 does not depend on the choice of $m$ and $j$. Since  $\alpha\mu>\gamma-1$ implies $\frac{\mu(\alpha +m)}{\gamma}>\frac{\gamma-1}{\gamma}$,   Remark \ref{integralomega} yields  that $\int_{\Omega}N^{k-\gamma}_\Omega |h^j_{k}|^2\omega_0^d$ is constant for all $j$ and thus (\ref{integral}) does not depend on $j$ (observe also that for $j=0$ one obtains the term $s_m$  of $s$ in Lemma \ref{immersion} and for $j=m=0$ we recover the first term of $f$, $f_0=1$). Thus if $\alpha\,g(\mu)$ is balanced the quantity
$$\frac{\pi}{(d+1)!} \frac{(m+ \alpha-1)!}{(\alpha-1)!}\frac{1}{(\alpha-(d+2)+1)\cdots(\alpha-(d+2)+m)}\int_{\Omega}N^{\mu(\alpha +m)-\gamma}_\Omega\omega_0^d$$
does not depend on $m$. By \cite[Prop. 2.1, p. 358]{roosformula} this is equivalent to ask that the quantity
$$\frac{(m+ \alpha-1)!}{(\alpha-d-1)\cdots(\alpha-d+m-2)}\frac{F(\mu(\alpha +m)-\gamma)}{F(0)}$$
does not depend on $m$, where 
$$\frac{F(s)}{F(0)}=)\prod_{j=1}^r\frac{\Gamma\left(s+1+\frac{(j-1)a}{2}\right)\Gamma\left(b+2+\frac{(r+j-2)a}{2}\right)}{\Gamma\left(1+\frac{(j-1)a}{2}\right)\Gamma\left(s+b+2+\frac{(r+j-2)a}{2}\right)} ,$$
for $\Gamma$ the usual Gamma function and $(a,b,r)$ the domain's invariants described in Table \ref{inv}. Deleting the terms which do not depends on $m$ and changing the orders of terms in the argument of the Gamma functions, we get 
$$(\alpha +m-1)\cdots(\alpha +m-d)\prod_{j=1}^r\frac{\Gamma\left(\mu(\alpha +m)-\gamma+1+\frac{(j-1)a}{2}\right)}{\Gamma\left(\mu(\alpha +m)-\gamma+b+2+\frac{(r+j-2)a}{2}\right)},$$
$$(\alpha +m-1)\cdots(\alpha +m-d)\prod_{j=1}^r\frac{\Gamma\left(\mu(\alpha +m)-\gamma+1+\frac{(j-1)a}{2}\right)}{\Gamma\left(\mu(\alpha +m)-\gamma+b+2+\frac{(2r-j-1)a}{2}\right)},$$
$$(\alpha +m-1)\cdots(\alpha +m-d)\prod_{j=1}^r\frac{\Gamma\left(\mu(\alpha +m)-\gamma+1-\frac{a}{2}+\frac{ja}{2}\right)}{\Gamma\left(\left[\mu(\alpha +m)-\gamma+1-\frac{a}{2}+\frac{ja}{2}\right]+b+1+ra-ja\right)}.$$
Since the quantity $b+1+ra-ja$ is a positive integer, by the well-known property $\Gamma(z+1)=z\Gamma(z)$ we get
\begin{equation}\label{final}
\frac{(\alpha +m-1)\cdot\ldots\cdot(\alpha +m-d)}{\underbrace{\left(A+b+a(r-1)\right)\cdot\ldots\cdot A}_{j=1,\ b+1+a(r-1) \,\textrm{terms}}\cdot\ldots\cdot\underbrace{\left(A+b\right)\cdot\ldots\cdot A}_{j=r,\ b+1\, \textrm{terms}}},
\end{equation}
where $A=\mu(\alpha +m)-\gamma+1-\frac{a}{2}+\frac{ja}{2}$. A necessary condition for the above quantity to be indipendent from $m$ is that numerator and denominator regarded as polynomials in $m$ have the same degree, i.e.
$$d=\frac{(b+1+a(r-1))!}{b!}.$$
Observe that such condition is satisfied by the complex hyperbolic space whose rank is $r=1$ and $b=d-1$. Althought no one of the other domains satisfies it. It remains to show that $\left(M_{\CH^d(\mu)}, \alpha \, g(\mu)\right)$ is  balanced if and only if $\mu=1$. By (\ref{final}), when $\Omega =\C H^{d}$, 
i.e.  $r=1$, $b=d-1$ and $\gamma=d+1$, we get
\begin{equation}
\frac{(\alpha +m-1)\cdot\ldots\cdot(\alpha +m-d)}{\left(\mu(\alpha +m)-1\right)\cdot\ldots\cdot \left(\mu(\alpha +m)-d\right)},\nonumber
\end{equation}
that is independent from $m$ if and only if $\mu=1$, as wished. 
\end{proof}

\end{document}